\documentclass[conference]{IEEEtran}
\pdfoutput=1
\usepackage{flushend}
\usepackage[T1]{fontenc}
\usepackage[utf8]{inputenc}
\usepackage{placeins}
\usepackage{amsfonts}
\usepackage[hidelinks]{hyperref}
\usepackage{makecell}
\usepackage{caption}
\ifCLASSOPTIONcompsoc
\usepackage[caption=false,font=footnotesize,labelfont=sf,textfont=sf]{subfig}
\else
\usepackage[caption=false,font=footnotesize]{subfig}
\fi

\usepackage{graphicx,amssymb,amsmath,array,url}

\newtheorem{definition}{Definition}

\newtheorem{theorem}{Theorem}
\newtheorem{corollary}{Corollary}

\newcommand{\twocases}[4]{%
  \begin{cases}
    #1 & #2\\
    #3 & #4
  \end{cases}%
}

\renewcommand{\epsilon}{\varepsilon}
\IEEEoverridecommandlockouts    % to create the author's affliation portion
\graphicspath{{./PDFs/}}
\DeclareGraphicsExtensions{.pdf}
\begin{document}
\title{\ \\ \LARGE\bf Encoding and Visualization in the Collatz Conjecture \thanks{George M. Georgiou is with the School of Computer Science and Engineering, California State University, San Bernardino, CA 92407, USA (email: georgiou@csusb.edu).}}
\author{George M. Georgiou}

\maketitle
\thispagestyle{plain}
\pagestyle{plain}

\begin{abstract}
The Collatz conjecture is one of the easiest mathematical problems to state and yet it remains unsolved. For each $n\ge 2$ the Collatz iteration is mapped to a binary sequence and a corresponding unique integer which can recreate the iteration. The binary sequence is used to produce the Collatz curve, a 2-D visualization of the iteration on a grid, which, besides the aesthetics, provides a qualitative way for comparing iterations. Two variants of the curves are explored, the r-curves and on-change-turn-right curves.  There is a scarcity of acyclic r-curves and only three r-curves were found having a cycle of minimum length greater than 4.
\end{abstract}
\section{Introduction}
\IEEEPARstart{T}{he} Collatz Conjecture states that starting with a positive integer $n$, the iteration of the function
\begin{equation}
  \label{eq:1}
  C(n) = \twocases{\frac{n}{2},}{\text{if $n\equiv 0 \pmod 2 $}}{3n+1,}{\text{if $n\equiv 1 \pmod 2$,}}
\end{equation}
eventually reaches 1; that is to say that there exists a finite integer $k$ such that $C^{k}(n) =1$. For example, starting with $n=3$, $C^0(3)=3$, we obtain the sequence $3 \rightarrow 10 \rightarrow 5 \rightarrow 16 \rightarrow 8 \rightarrow 4 \rightarrow 2 \rightarrow 1$. Continuing the iteration, it is noted that the sequence is caught in the cycle $1 \rightarrow 4 \rightarrow 2 \rightarrow 1$. The Collatz conjecture can be stated in another way: starting the iteration with any  positive integer $n$  there is no other cycle than this one and the iteration does not diverge to infinity.

The origins of the Collatz Conjecture can be traced back to the 1930's to the German mathematician Luther Collatz \cite{Lagarias1985}. Although known in the 1950's among mathematicians, it first appeared in print in 1971 \cite{Coxeter1971,Lagarias1985}. Despite the simplicity with which it can be stated and despite the multitude of mathematicians who studied it and attempted to solve it, its proof has remained elusive. Empirically, it has been verified to be true for $n < 5 \times 2^{60} \approx 5.7646 \times 10^{18}$ \cite{OliveiraeSilva2010}.

The consensus seems to be that a proof is not within reach at this time. In general, randomness and unpredictability govern the sequences produced by the iterations from one number to the next. Paul  Erdős famously commented ``Mathematics may not be ready for such problems'' \cite{Lagarias1985}.  The conjecture, as it circulated in many campuses through many mathematicians, it came to be known by various names such as the \emph{$3x+1$ problem,} the \emph{Syracuse Problem}, \emph{Hasse's Algorithm}, \emph{Kakutani's Problem} and \emph{Ulam's Problem} \cite{Lagarias1985,Lagarias2006}. The sequence of numbers produced by the Collatz iteration is known as  \emph{hailstone numbers} \cite{Hayes1984}.

It is noted that in Equation~\ref{eq:1}, when the second rule is applied, the resulting number $3n+1$ is necessarily an even number, and hence the iteration can be written more compactly thusly \cite{Terras1976,Terras1979,Lagarias1985}:
\begin{equation}
  \label{eq:2}
  T(n) =
  \begin{cases}
    \frac{n}{2}, & {\text{if $n\equiv 0 \pmod 2 $}}\\
    \frac{3n+1}{2}, & {\text{if $n\equiv 1 \pmod 2$.}}
  \end{cases}
\end{equation}
This version of the iteration is widely used and it will be used for the rest of this paper. Using the same example $n=3$, with iteration $T^k(n)$ the sequence produced is $3 \rightarrow 5  \rightarrow 8 \rightarrow 4 \rightarrow 2 \rightarrow 1$, and the cycle becomes $1 \rightarrow 2 \rightarrow 1$.

The \emph{total stopping time} $\sigma_\infty(n)$ of $n$ is the least $k$ applications of $T$ that the sequence of the iteration will reach 1 for the first time:
\begin{equation}
  \label{eq:3}
  \sigma_\infty(n) = \inf\:\{k:T^k(n) =1\}.
\end{equation}
For example, $\sigma_\infty(3)= 5.$

In this paper, the iteration $T^k(n)$ of each number $n$ is encoded as a binary number (or string), and in turn as a unique positive integer $m$. The binary string is mapped to a curve in 2-D in an analogous way that the binary string of the Dragon Curve is mapped to the corresponding curve.  
\section{The Binary Encoding}
\label{sec:binary-string-}

The conversion of a decimal number $n$ to its corresponding binary representation involves successive divisions by 2 and looking at the parity of the result.  We note that the process is similar to the Collatz iteration: each time $T(n)$ is applied a decision is made based on whether $n$ is even or odd. We choose to have a ``1'' if $T^k(n)$ is even and ``0'' if odd. For example, $n=3$ is encoded as 11100.  The bits are produced, just like in the usual decimal-to-binary conversion, from right to left. Table~\ref{tab:1} shows the binary encoding and decimal equivalent.

For all $n \geq 2$, assuming that the Collatz Conjecture is true, the last application of $T$ is on $2$, an even number, before reaching $1$. Making the same assumption, the binary encoding of all $n \geq 5$ begins with $111$, noting that  eventually  the iteration will end with $8 \rightarrow 4 \rightarrow 2 \rightarrow 1$. The length of the binary encoding of $n$ is $\sigma_\infty(n)$ since the number of digits is exactly the number of applications of $T$ needed to reach $1$ in the sequence of iterations. The choice of when to use  ``1'' or ``0'' was made so that the binary encoding of $n$, the original number, starts with a ``1'', to avoid having any leading zeros. Had we had any leading zeros, the binary string would not have been recoverable from its corresponding decimal number, since leading zeros are not significant and are ignored when converting to decimal. This sequence was independently published as sequence A304715 in OEIS \cite{A304715}, which we noticed after a draft of this paper had been released.

Given the binary encoding of $n$, it is a simple matter to find the original $n$. It can be done by the following process: Scan its binary encoding $B=b_lb_{l-1}\ldots{}b_1b_0$, where $b_l$ is the most significant bit, from left to right. Starting with initial $n=1$ and $i=0$, iteratively apply $R^i(n,b)$, for $i = 0, 1, \ldots, l$:
\begin{equation}
  \label{eq:5}
  R^i(n,b) =
  \begin{cases}
    2n, & {\text{if $b_{l-i} = 1$}}\\
    \frac{2n-1}{3}, & {\text{if $b_{l-i}=0$}}.
  \end{cases}
\end{equation}
As an example, given the encoding $B=11100$, starting with $n=1$ and iteratively applying $R^i(n,b)$ we obtain the sequence $1 \rightarrow 2 \rightarrow 4 \rightarrow 8 \rightarrow 5 \rightarrow 3$, and thus we obtain $n=3$. This is the reverse sequence produced by Equation~\ref{eq:2}.

\begin{table}[ht]
  \centering\small
  \begin{tabular}{|r|r|r|}\hline
$n$&\makecell[cc]{Binary encoding}&\makecell[cc]{Decimal}\\\hline
2 & 1 & 1 \\
3 & 11100 & 28 \\
4 & 11 & 3 \\
5 & 1110 & 14 \\
6 & 111001 & 57 \\
7 & 11101101000 & 1896 \\
8 & 111 & 7 \\
9 & 1110110100010 & 7586 \\
10 & 11101 & 29 \\
11 & 1110110100 & 948 \\
12 & 1110011 & 115 \\
13 & 1110110 & 118 \\
14 & 111011010001 & 3793 \\
15 & 111011110000 & 3824 \\
16 & 1111 & 15 \\
17 & 111011010 & 474 \\
18 & 11101101000101 & 15173 \\
19 & 11101101001100 & 15180 \\
20 & 111011 & 59 \\
21 & 111110 & 62 \\
22 & 11101101001 & 1897 \\
23 & 11101111000 & 1912 \\
24 & 11100111 & 231 \\
25 & 1110110100110010 & 60722 \\
26 & 11101101 & 237 \\
27 & \makecell[tl]{%
     111011110001101110\\
    101011100001100000\\
    010010001101000010\\
    \hphantom{01}0010010100000100} & 1102691417057682138372 \\
28 & 1110110100011 & 7587 \\
29 & 1110110100110 & 7590 \\
    30 & 1110111100001 & 7649\\\hline
  \end{tabular}
  \caption{The binary and corresponding decimal encoding.}\label{tab:1}
\end{table}

\section{Collatz Curves: Visualizing the Binary Encoding}
\label{sec:visualizing-n}
For each positive integer $n\ge 2$, we will visualize the Collatz sequence generated by the iteration $T^k(n)$ by making use of its binary encoding (Table~\ref{tab:1}).  The binary encoding is viewed as a random walk where each of its digits generates a unit length line segment in a certain direction. 
 
\begin{figure}[!htbp]
  \centering
  \includegraphics[width=0.5\linewidth,keepaspectratio]{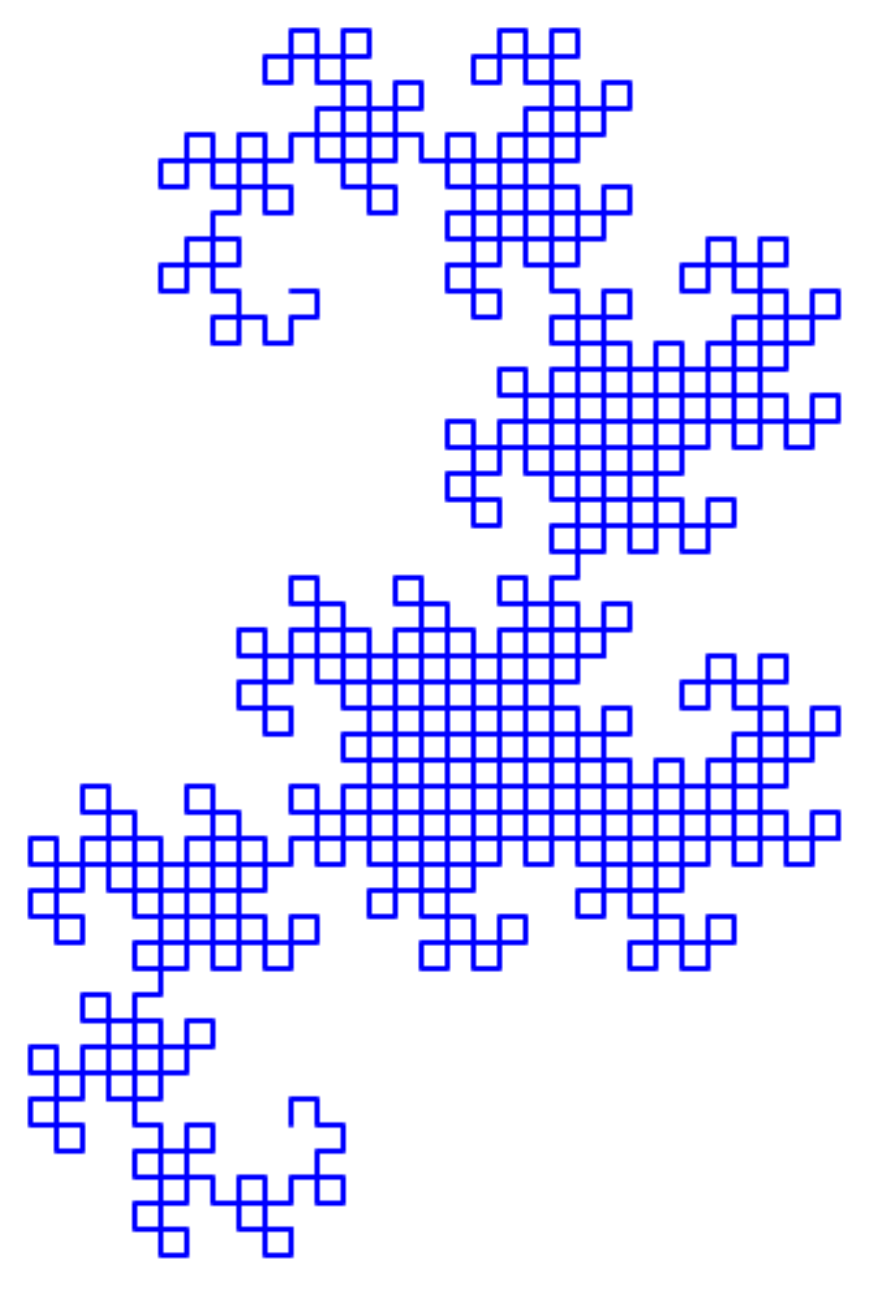}
  \caption{The Dragon Curve of order $10$.}
  \label{fig:1}
\end{figure}

The manner of generating the Collatz curve corresponding to $n$ is similar to that of generating the Dragon Curve (Figure~\ref{fig:1})  \cite{Davis1970,Gardner1978a}.
The string $S_k$ of the Dragon Curve of order $k$ is generated as follows. There are two symbols, R (right) and L (left). In the first iteration (order $k=1$) the string $S_1$ is just the symbol R. In each successive iteration the new string $S_{n+1}$ is recursively produced:
\begin{equation}
S_{n+1}=S_n + R + \overline{S}_n,
\end{equation}
where $+$ is string concatenation and $\overline{S}_n$ is $S_n$ reversed with each of its symbols flipped.  Table~\ref{tab:2} shows the first four iterations.

\begin{table}[ht]
  \centering\small
  \begin{tabular}{|c|c|}\hline
Iteration ($k$)&\makecell[cc]{String $S_k$}\\\hline
1&R\\
2&R R L\\
3&R R L R R L L\\
4&R R L R R L L R R R L L R L L\\\hline
  \end{tabular}
  \caption{The first four iterations of the Dragon Curve.}\label{tab:2}
\end{table}

Let $S_k=s_ls_{l-1}\ldots s_1s_0.$ From $S_k$ a curve  can be produced such as the one in Figure~\ref{fig:1} ($k=10$) by making use of turtle graphics \cite{Abelson1986,dragon}: string $S_k$ is scanned from left to right:\\ 
\begin{center}
\begin{tabular}[h]{ll}
  if $s_i = R$, &turn right $90^\circ$, move forward one unit\\
  if $s_i=L$,& turn left $90^\circ$, move forward one unit.  
\end{tabular}
\end{center}

   The initial orientation of the turtle in all figures in this paper is facing east.

The curve of each Collatz sequence is produced in an analogous way. The digits of the binary encoding $B$ of $n$ are scanned from right to left.
\begin{center}
\begin{tabular}[h]{ll}
  if $b_i = 1$, &turn right $90^\circ$, move forward one unit\\
  if $b_i=0$,& turn left $90^\circ$, move forward one unit.  
\end{tabular}
\end{center}
 
\begin{figure}[!htbp]
  \centering
  \includegraphics[height=1.0cm,keepaspectratio]{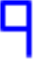}
  \caption{The curve for $n=3$.}
  \label{fig:2}
\end{figure}

\begin{figure}[!htbp]
  \centering
  \includegraphics[width=0.3\linewidth,keepaspectratio]{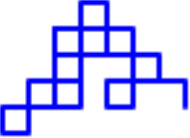}
  \caption{The curve for $n=27$.}
  \label{fig:3}
\end{figure}

\begin{figure}[!htbp]
  \centering
  \includegraphics[width=0.5\linewidth,keepaspectratio]{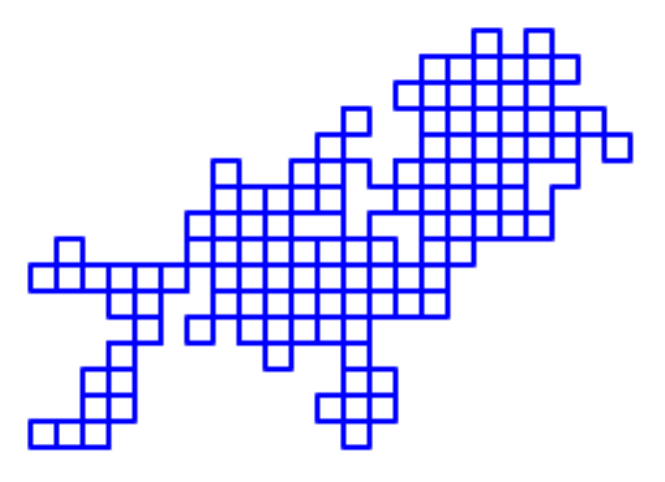}
  \caption{The ``lion curve'' for $n=75,128,138,247$.}
  \label{fig:4}
\end{figure}

The curve for $n=3$ appears in Figure~\ref{fig:2}. Its binary encoding is $11100$. The initial point, that is the initial position of the turtle, is at the bottom of the curve. A more involved curve is generated for $n=27$ (Figure~\ref{fig:3}), a number known for its relatively long binary encoding, which is of length $70$. A still more involved curve, for $n=75,128,138,247$, is shown in Figure~\ref{fig:4}, which bears a resemblance of the figure of a lion. Its binary sequence is of length $767$.

\begin{figure}[!htbp]
  \centering
  \includegraphics[width=0.10\linewidth,keepaspectratio]{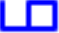}
  \caption{The curve for $n=6$.}
  \label{fig:5}
\end{figure}

\begin{figure*}[ht]
  \centering
  \includegraphics[width=0.95\textwidth,keepaspectratio]{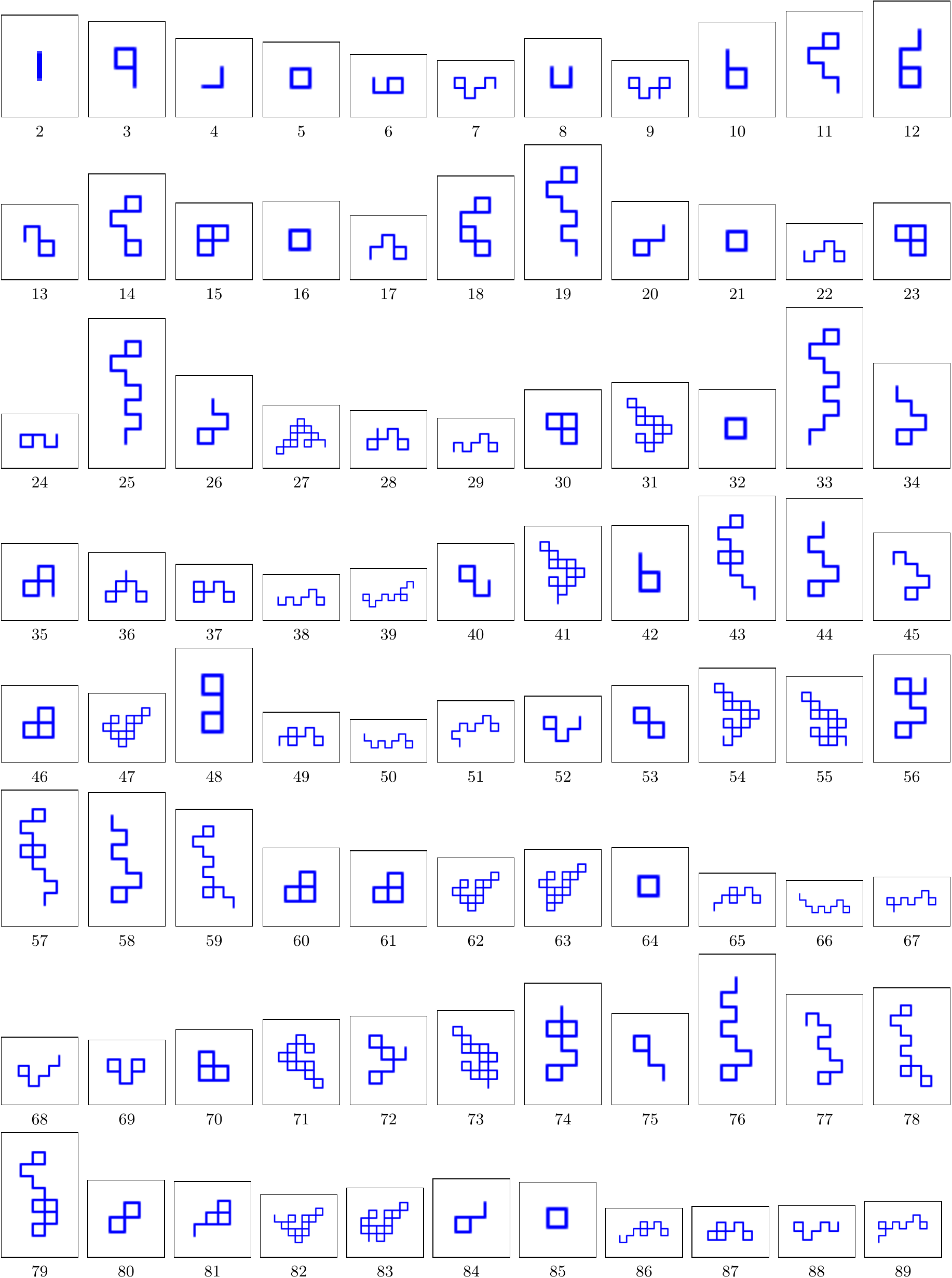}
  \caption{The Collatz curves for $n=2, 3, \ldots, 89.$}
  \label{fig:6}
\end{figure*}

The curve for $n=6$ (Figure~\ref{fig:5}, also in Figure~\ref{fig:6}), whose binary encoding is $111001$, coincidentally forms the shape of $6$! The orientation is arbitrary since it only depends on the starting orientation of the turtle.

\begin{figure*}[ht]
  \centering
  \includegraphics[width=0.95\textwidth,keepaspectratio]{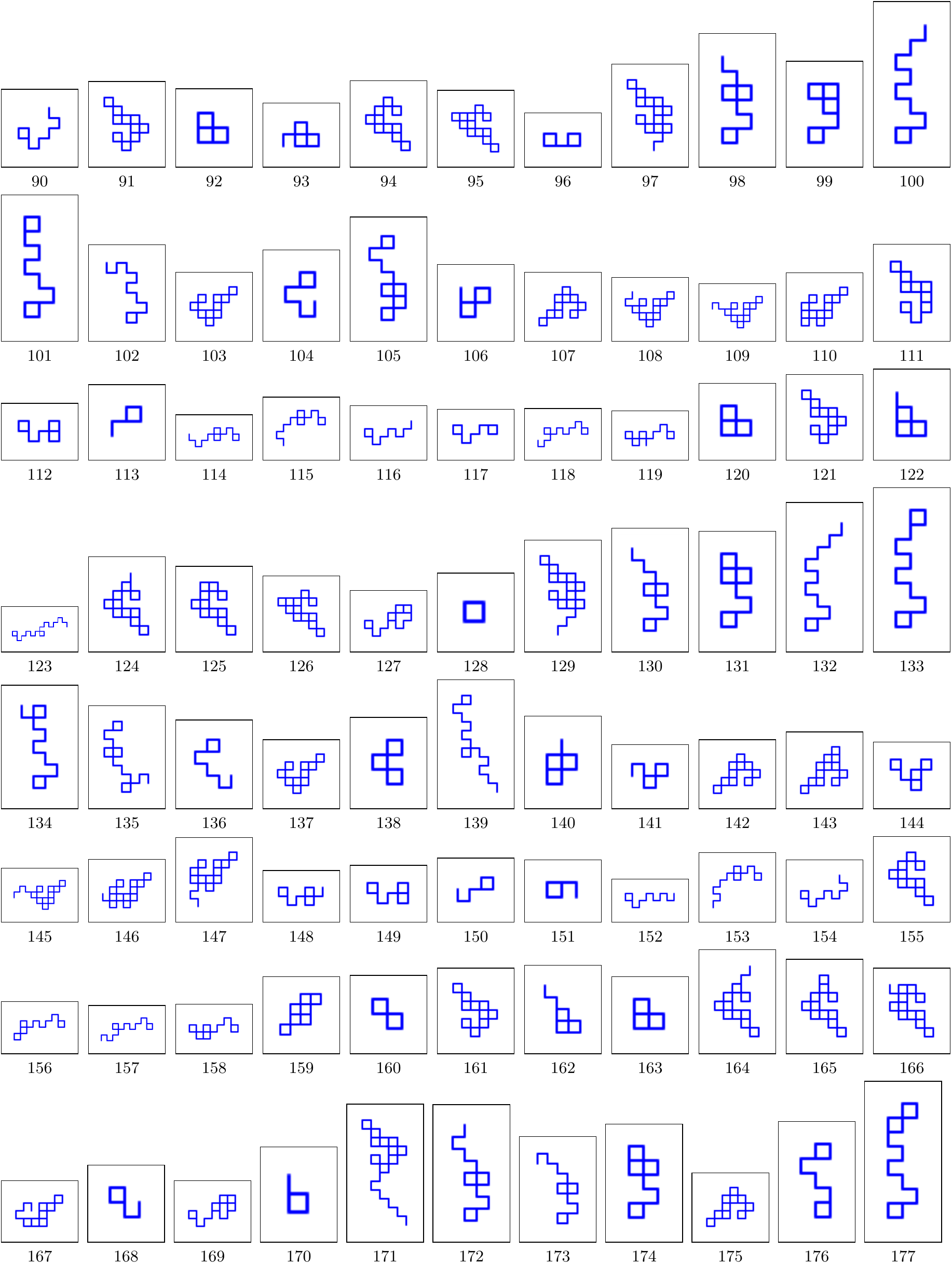}
  \caption{The Collatz curves for $n=90, 91, \ldots, 177.$}
  \label{fig:7}
\end{figure*}

Figure~\ref{fig:6} shows the Collatz curves for $2\le n \le 89$ and Figure~\ref{fig:7} those for $90\le n \le 177$.

One important difference between the Dragon Curve and the curves we generate using the Collatz iteration is the former does not intersect itself \cite{Davis1970}, whereas the latter ones do.  It is possible, and it often happens, that the curve as it is generated  writes over previously generated segments. Hence, although the binary encoding is unique for each $n$, as discussed in Section~\ref{sec:binary-string-}, its corresponding curve is not. As an example, the curves for $n=2^k, k \ge 4$ appear identical: a unit square, which can be seen in  Figures~\ref{fig:6} and \ref{fig:7} for $n=16, 32, 64, 128.$ But powers of two are not the only ones that produce a square, e.g. consider the figures of $n= 5, 21, 85.$

\begin{theorem}\label{theo:1}
  Assuming that the Collatz Conjecture is true, all curves for $n \ge 2$, except those for $n = 2,4, \text{or } 8$, end with a cycle, i.e. a unit square.
\end{theorem}
\begin{IEEEproof}
For $n=2,4,8$ the binary encoding is of length three or less (Table~\ref{tab:1}), which precludes forming a cycle, since a cycle requires at least 4 segments. For all other $n$ the binary encoding is of length 4 or more and the iteration ends in the cycle $4\rightarrow 2\rightarrow 1$. Necessarily, the sequence before reaching $1$ for the first time is $8\rightarrow 4\rightarrow 2\rightarrow 1$, which corresponds to $111$.  Right before $111$ is drawn, the turtle is at the end of an existing segment. Turning right and drawing 3 times creates a final cycle.
\end{IEEEproof}
\begin{definition}
  The \emph{girth} $g(G)$ of a graph $G$ is the minimum length of a cycle contained in $G$ \cite{Diestel2017}.  If $G$ has no cycles, $g(G)=\infty$.
\end{definition}

\begin{corollary}
  When the Collatz curve is viewed as a graph $G(n)=(V,E)$, $g(G(n))=4$ for all $n > 8$.
\end{corollary}

The curves allow for easy identification of similar patterns in the Collatz iteration. As an example, visually inspecting the curves in Figure~\ref{fig:6}, we note that those for $n=27$, $31$, $41$, $47$, $54$, $55$, $62$, $63$, $71$, $73$, $82$, and $83$ bear a strong resemblance to each other. A closer look at the binary encodings reveals that indeed when compared pairwise they have a long common substring (subsequence). Table~\ref{tab:4} shows the lengths of the longest common substrings, all of which happen to start from the leftmost bit position (see Table~\ref{tab:5}). Each common substring is of length at least 60. In Table~\ref{tab:4}, on the diagonal is the length of the binary encoding of the corresponding $n$.
\begin{table}[ht]
  \centering%\scriptsize
\resizebox{\columnwidth}{!}{\begin{tabular}{c|cccccccccccc}
    &27&31&41&47&54&55&62&63&71&73&82&83\\\hline
  27&\textbf{70}&67&69&66&70&66&68&60&65&66&69&66\\
  31&67&\textbf{67}&67&66&67&66&67&60&65&66&67&66\\
  41&69&67&\textbf{69}&66&69&66&68&60&65&66&69&66\\
  47&66&66&66&\textbf{66}&66&66&66&60&65&66&66&66\\
  54&70&67&69&66&\textbf{71}&66&68&60&65&66&69&66\\
  55&66&66&66&66&66&\textbf{71}&66&60&65&71&66&70\\
  62&68&67&68&66&68&66&\textbf{68}&60&65&66&68&66\\
  63&60&60&60&60&60&60&60&\textbf{68}&60&60&60&60\\
  71&65&65&65&65&65&65&65&60&\textbf{65}&65&65&65\\
  73&66&66&66&66&66&71&66&60&65&\textbf{73}&66&70\\
  82&69&67&69&66&69&66&68&60&65&66&\textbf{70}&66\\
  83&66&66&66&66&66&70&66&60&65&70&66&\textbf{70}
\end{tabular}}
\caption{The pairwise length of the longest common substring of the binary encodings. On the diagonal is the length of the binary encoding of the corresponding $n$.}\label{tab:4}
\end{table}

\begin{table*}[ht]
  \centering\normalsize%small
  % \resizebox{\columnwidth}{!}{
  \begin{tabular}{c|l}
    n&\makecell[cc]{Binary encoding}\\\hline
27&11101111000110111010101110000110000001001000110100001000100\textbf{1}0100000100\\
31&11101111000110111010101110000110000001001000110100001000100\textbf{1}0100000\\
41&11101111000110111010101110000110000001001000110100001000100\textbf{1}010000010\\
47&11101111000110111010101110000110000001001000110100001000100\textbf{1}010000\\
54&11101111000110111010101110000110000001001000110100001000100\textbf{1}01000001001\\
55&11101111000110111010101110000110000001001000110100001000100\textbf{1}01000011000\\
62&11101111000110111010101110000110000001001000110100001000100\textbf{1}01000001\\
63&11101111000110111010101110000110000001001000110100001000100\textbf{1}11000000\\
71&11101111000110111010101110000110000001001000110100001000100\textbf{1}01000\\
73&11101111000110111010101110000110000001001000110100001000100\textbf{1}0100001100010\\
82&11101111000110111010101110000110000001001000110100001000100\textbf{1}0100000101\\
83&11101111000110111010101110000110000001001000110100001000100\textbf{1}0100001100\\
  \end{tabular}%}
\caption{Numbers with similar binary encodings, i.e. similar trajectories in the Collatz iteration. They have been identified visually via their Collatz curves. Pairwise they have a common substring, leftmost bits, of length at least 60. The 60th bit from the left is in bold.}\label{tab:5}
\end{table*}
 
\section{Reverse Curves}
\label{sec:reverse-curves}

\begin{figure}[!htbp]
  \centering
  \fbox{\includegraphics[height=2.0cm,keepaspectratio]{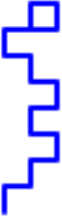}}\hspace*{2.0cm}
    \fbox{\includegraphics[height=2.0cm,keepaspectratio]{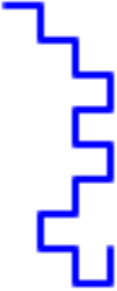}}
  \caption{For $n=33$, the Collatz curve (left) and the r-curve (right). Note the absence of a cycle in the r-curve.}
  \label{fig:8}
\end{figure}

\begin{figure*}[ht]
  \centering
  \includegraphics[width=0.9\textwidth,keepaspectratio]{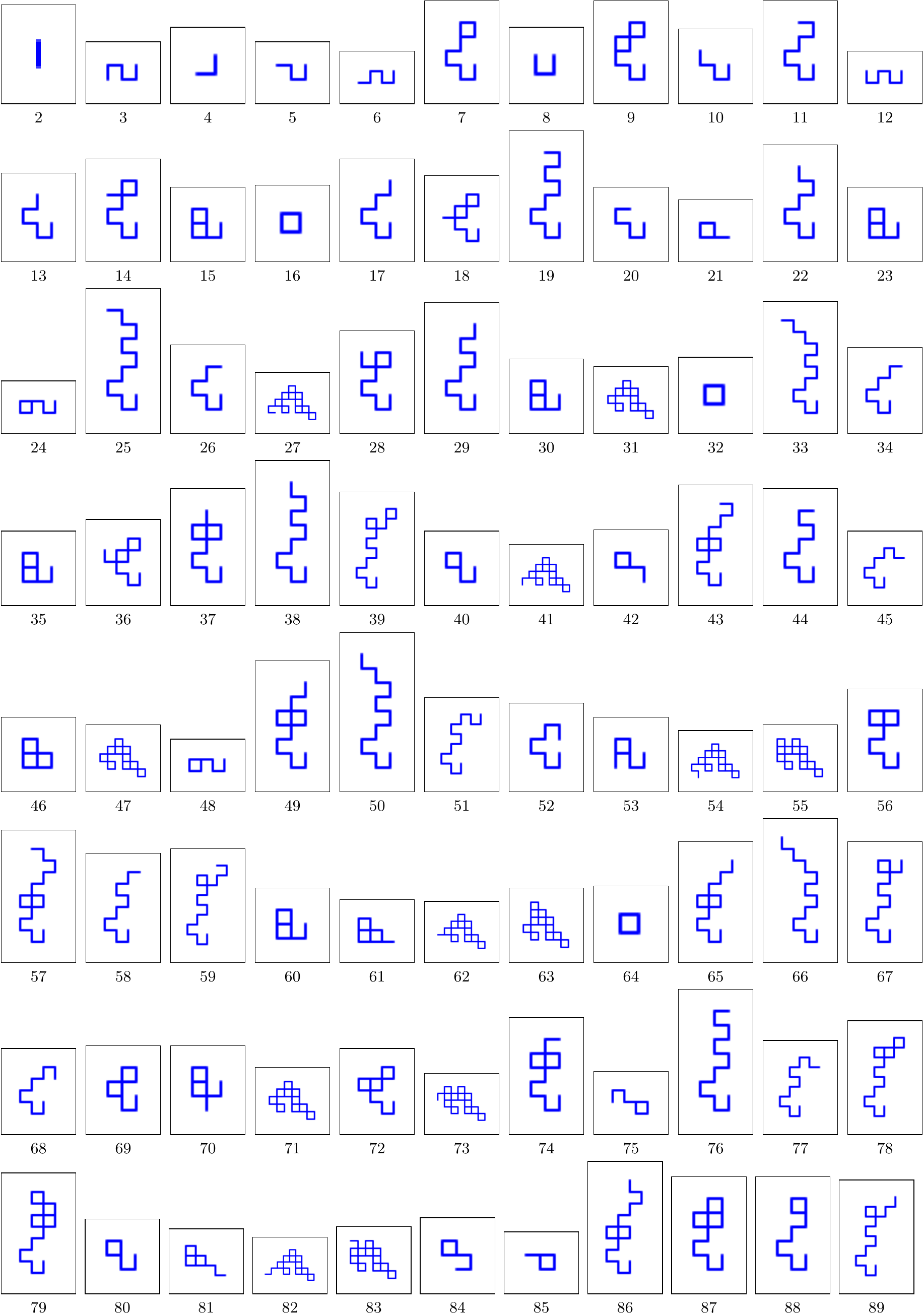}
  \caption{The r-curves for $n=2, 3, \ldots, 89.$}
\label{fig:13}
\end{figure*}

\begin{figure*}[ht]
  \centering
  \includegraphics[width=0.9\textwidth,keepaspectratio]{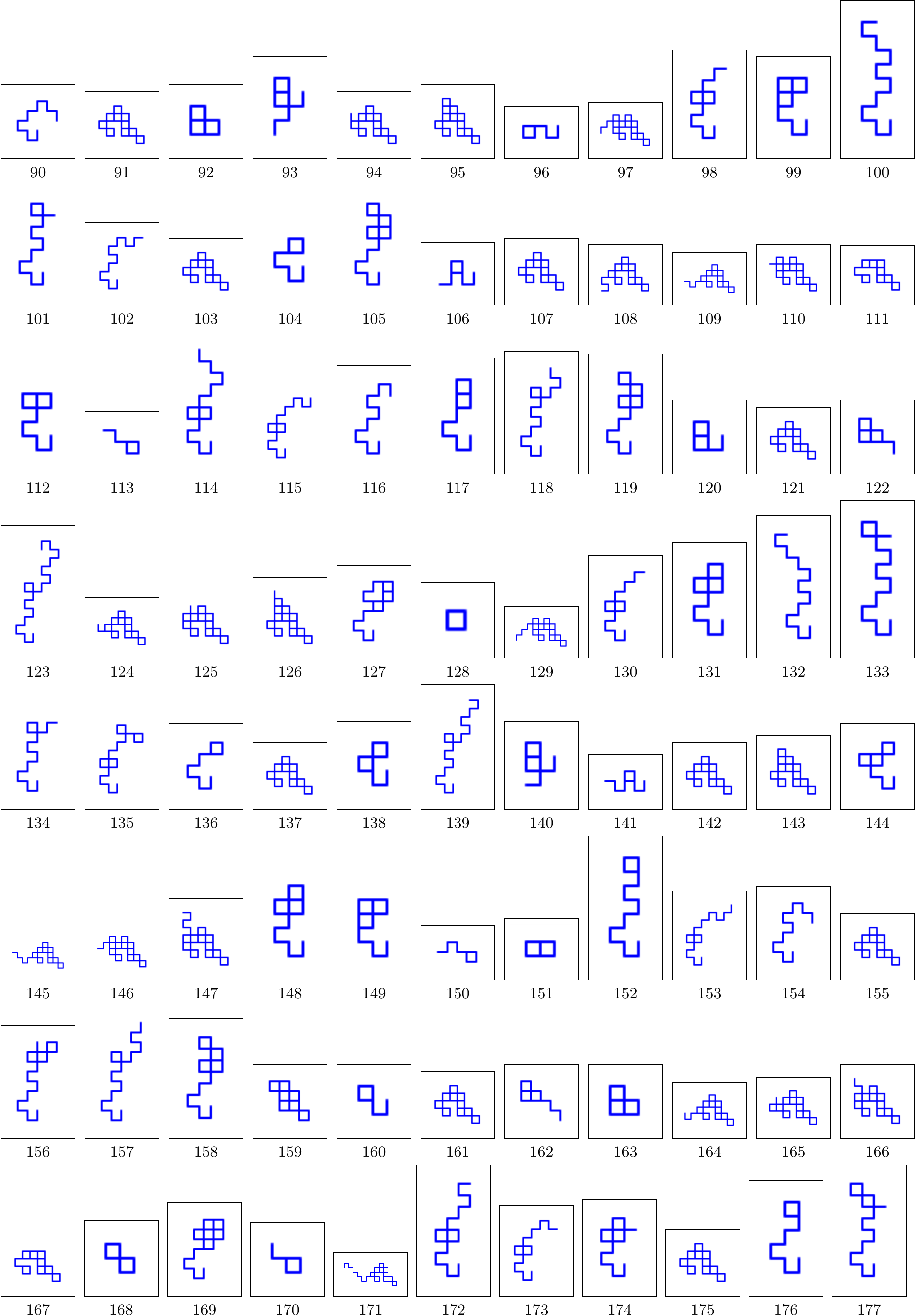}
  \caption{The r-curves for $n=90, 91, \ldots, 177.$}
\label{fig:14}
\end{figure*}

Another version of the Collatz curves can be obtained by applying the same rules except that the bits in the encoding of $n$ are scanned in reverse, from left-to-right. The resulting reverse curves, or \emph{r-curves}, are very similar to the Collatz curves.  In general they are not affine transformations of each other. Theorem~\ref{theo:1} does not apply, and thus the curves no longer necessarily have a cycle for $n>8$. As an example, compare the two curves for $n=33$ in Figure~\ref{fig:8}. The r-curve has no cycle. Figure~\ref{fig:13} shows the r-curves for $2\le n \le 89$ and Figure~\ref{fig:14} those for $90\le n \le 177$.

\begin{figure*}[!hbt]
  \centering
  \includegraphics[width=0.8\textwidth,keepaspectratio]{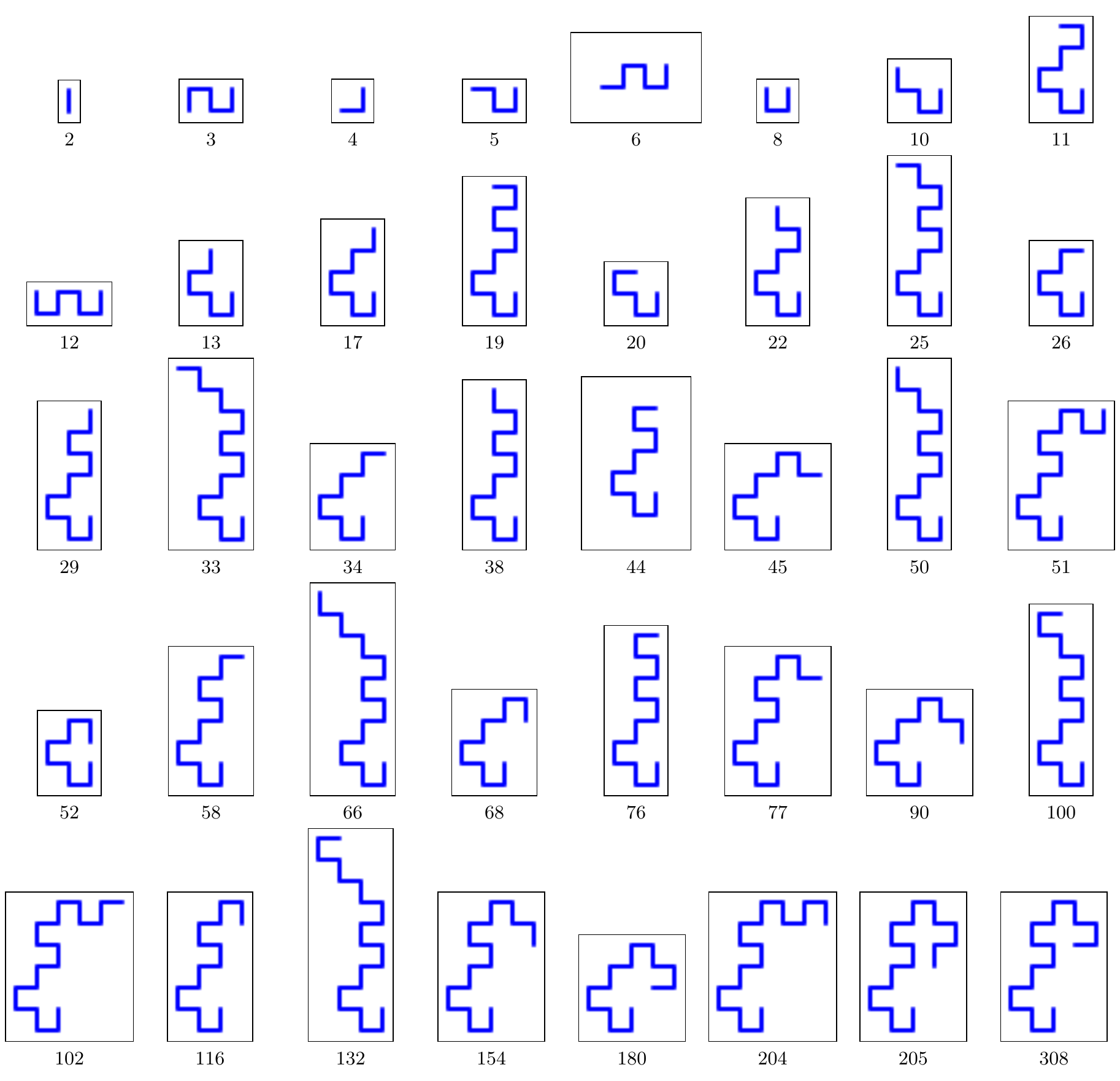}
  \caption{The r-curves without cycles.  Verified up to $n \le 10^8$.}\label{fig:9}
\end{figure*}

Two questions arise: (1) how many r-curves are devoid of cycles, and (2) how many r-curves have finite girth $g(G(n))>4$, i.e. they have a cycle of minimum length other than the unit square.  Using Python 3 and the software library \emph{NetworkX} \cite{Hagberg2008}, we have checked all r-curves up to $n \le 10^8$ and we found that there are only 40 r-curves with no cycles, with the largest $n=308$. All such curves can be found in Figure~\ref{fig:9}.

\begin{figure}[!tbhp]
\captionsetup[subfigure]{labelformat=empty}
\centering
\subfloat[273]{\fbox{\includegraphics[height=2.5cm]{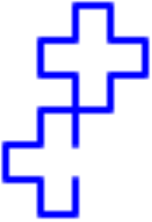}}}
\hfil
\subfloat[410]{\fbox{\includegraphics[height=2.5cm]{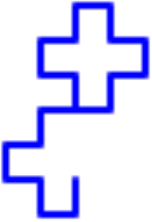}}}
\hfil
\subfloat[820]{\fbox{\includegraphics[height=2.5cm]{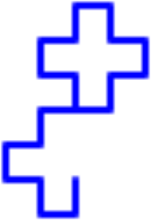}}}
\caption{Miracle curves: The only three r-curves found having a cycle of minimum length greater than $4$, verified up to  $n \le 10^8$. Their only cycle, the cross pattern, has length $12$. The curves for $410$ and $820$ are identical due to overwriting.}
\label{fig:10}
\end{figure}

In the same range, i.e. up to $n \le 10^8$, it was found that there are only three r-curves having finite girth $g(G(n))>4$, and those are for $n=273, 410, \text{and } 820$. In all three cases  $g(G(n))=12$. The only cycle in these curves is the cross pattern which can be seen in Figure~\ref{fig:10}. We call these curves \emph{miracle curves}. The ones for $410$ and $820$ are identical due to overwriting.
\section{Alternative Curves}
\label{sec:alternative-curves}

It is possible to use alternative schemes to draw the curves that correspond to the binary encoding of $n$ by modifying the rules. For example, we can use this rule for the turtle:\\

\begin{tabular}[h]{ll}
  if $b_i = b_{i-1}$,&  move forward one unit\\
  if $b_i \ne b_{i-1}$,& turn right $90^\circ$, move forward one unit.\\
\end{tabular}\\

\begin{figure}[!htbp]
\begin{center}
  \includegraphics[width=0.18\linewidth,keepaspectratio]{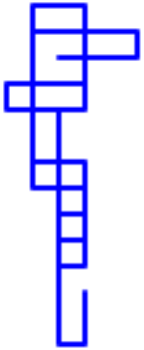}
  \captionof{figure}{The on-change-turn-right curve for $n=27$.}
  \label{fig:11}
\end{center}
\end{figure}

\begin{figure}[!htbp]
  \centering
  \includegraphics[width=0.3\linewidth,keepaspectratio]{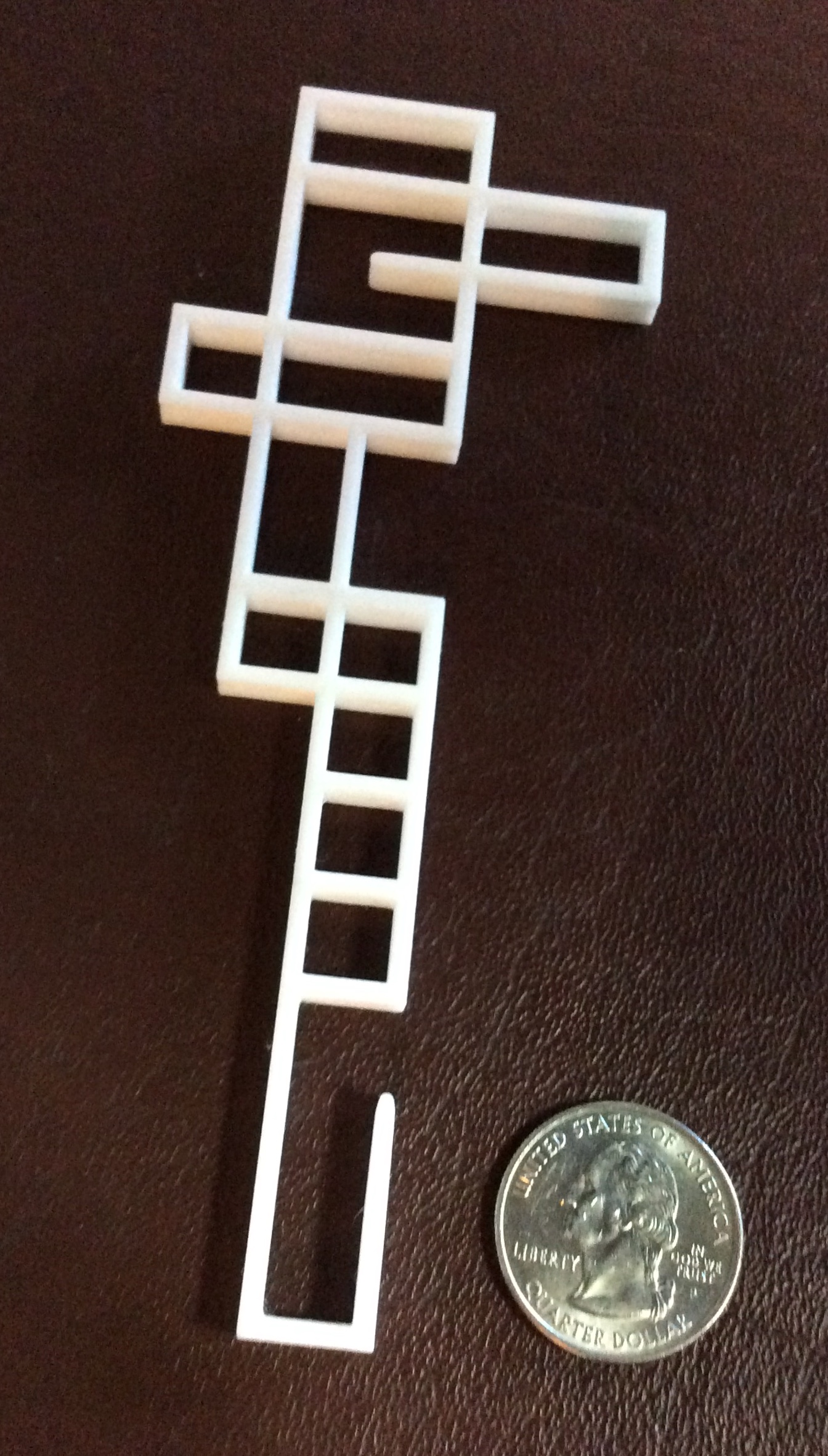}
  \caption{The on-change-turn-right curve for ${n=27}$ 3-D printed.}
  \label{fig:12}
\end{figure}

\begin{figure*}[ht]
  \centering
  \includegraphics[width=0.95\textwidth,keepaspectratio]{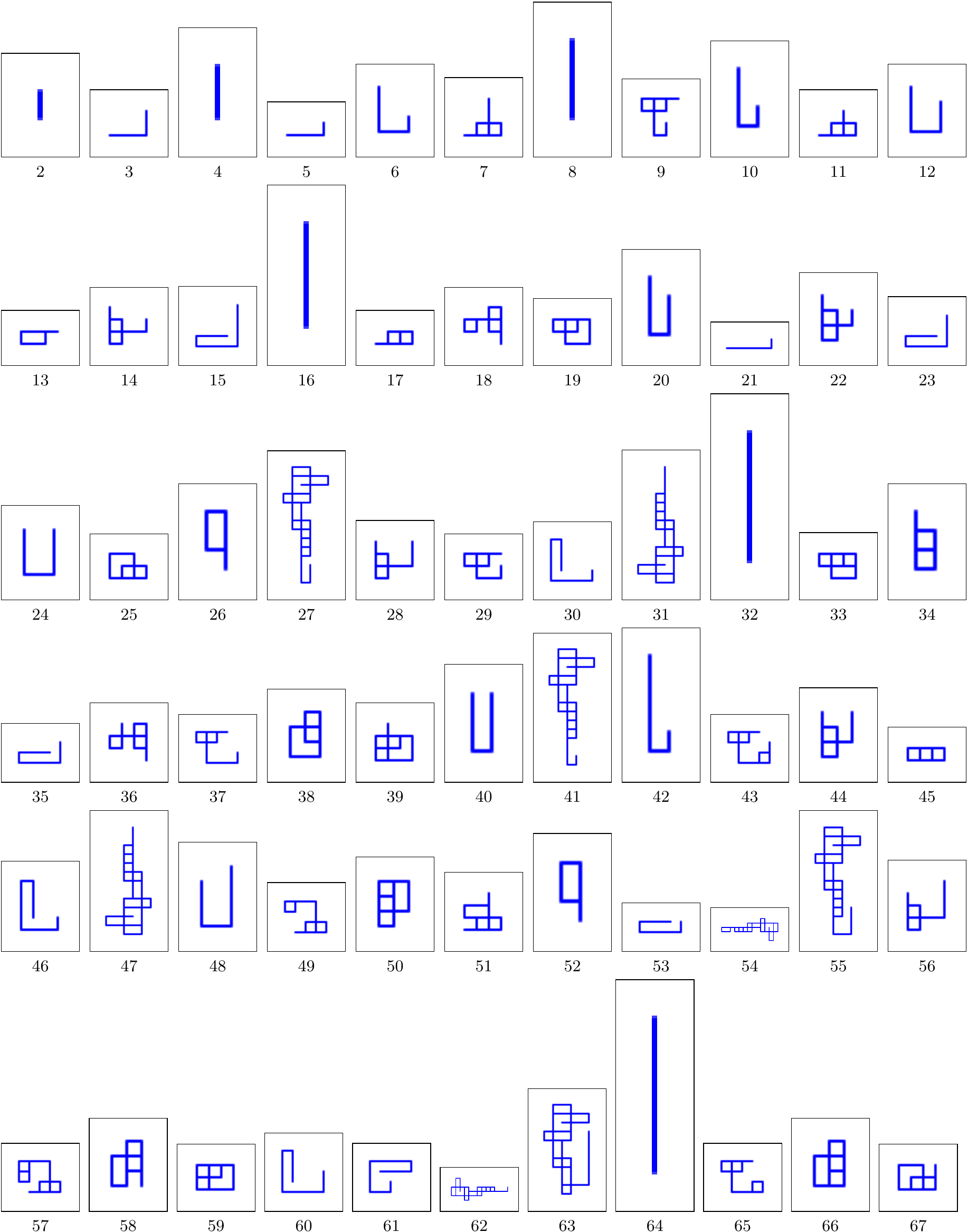}
  \caption{The on-change-turn-right curves for $n=2, 3, \ldots, 67.$}
  \label{fig:15}
\end{figure*}

\begin{figure*}[ht]
  \centering
  \includegraphics[width=0.95\textwidth,keepaspectratio]{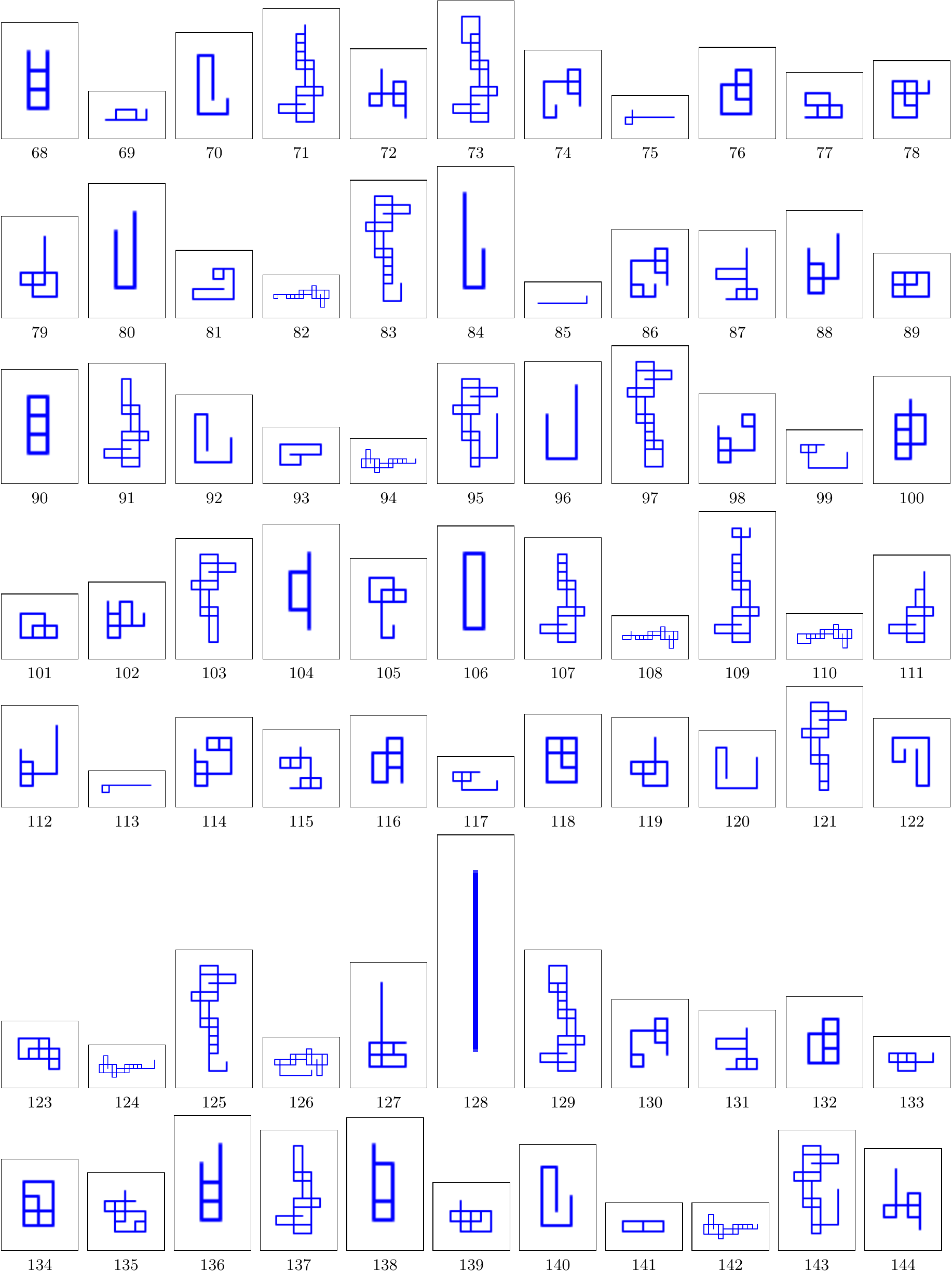}
  \caption{The on-change-turn-right curves for $n=68, 69, \ldots, 144.$}
  \label{fig:16}
\end{figure*}

\noindent In other words, as long as there is no change from one bit to the next, scanning the bits from right to left, either from 0 to 1 or vice versa, the turtle keeps moving forward. If there is a change, the turtle turns right $90^\circ$ and moves forward. The \emph{on-change-turn-right} curves produced with this rule appear in Figures~\ref{fig:15} and \ref{fig:16}. Visual identification of similar curves, e.g. those similar to the curve of $n=27$, is easily done. These curves are characterized by longer straight lines compared to the Collatz curves, Figures~\ref{fig:6} and \ref{fig:7}. For comparison with Figure~\ref{fig:3}, the on-change-turn-right curve for $n=27$ is shown in Figure~\ref{fig:11}. For purely aesthetic reasons, we show the picture of the same curve 3-D printed in Figure~\ref{fig:12}.


\section{Conclusion}
The Collatz iteration for each $n\ge 2$ is mapped to a binary sequence, which in turn is mapped to an integer. The binary sequence is translated to a 2-D curve on a grid, the Collatz curve, in an analogous way that the Dragon Curve is drawn from a binary sequence. Two variants of the curves were explored: the reverse curves (\emph{r-curves}) and the the \emph{on-change-turn-right} curves. There is a scarcity of acyclic r-curves, only 40 were found, and only three r-curves were found having a shortest length cycle of length greater than four, the \emph{miracle curves}. This was verified up to $n\le 10^8$. The Collatz curves, the r-curves, and the on-change-turn-right curves provide a way to visualize and qualitatively compare Collatz iterations.

\section*{Acknowledgment}
\label{sec:acknowledgement}
I would like to thank Rémy Sigrist for suggesting the exploration of the reverse curves and Arusyak Hovhannesyan for 3-D printing the on-change-turn-right curve for $n=27$.



\FloatBarrier

\renewcommand*{\UrlFont}{\rmfamily}
\bibliography{bib/collatz}
\end{document}